\documentclass[a4paper,11pt]{article}
\usepackage{amsmath, amssymb, mathtools, graphicx, parskip, amsthm, fancyhdr, dsfont, algorithm, algpseudocode, hyperref, bm, xcolor}

\newtheorem{thm}{Theorem}

\newtheorem{lemma}{Lemma}
\newtheorem{cor}{Corollary}
\theoremstyle{definition}

\newtheorem{remark}{Remark}

\newcommand{\bbP}{\mathbb{P}}
\newcommand{\bbE}{\mathbb{E}}
\newcommand{\bbN}{\mathbb{N}}
\newcommand{\bbR}{\mathbb{R}}

\title{A debiased Bernoulli factory and unbiased estimation of a probability}
\author{Jere Koskela \\
	\texttt{jere.koskela@tugraz.at} \\
	\small Institute of Statistics, \\
	\small Graz University of Technology and \\
	\small School of Mathematics, Statistics and Physics,\\
	\small Newcastle University
	\and
	Toni Karvonen\\
	\texttt{toni.karvonen@lut.fi}\\
	\small School of Engineering Sciences, \\
	\small Lappeenranta--Lahti University \\
	\small of Techonology LUT \\
	\and
	Krzysztof {\L}atuszy\'nski\\
	\texttt{k.g.latuszynski@warwick.ac.uk}\\
	\small Department of Statistics,\\
	\small University of Warwick
	\and
	Dario Span\`{o} \\
	\texttt{d.spano@warwick.ac.uk}\\
	\small Department of Statistics,\\
	\small University of Warwick
}
\date{\today}

\begin{document}

\maketitle

\begin{abstract}
Given a known function $f : [0, 1] \to (0, 1)$ and a random but almost surely finite number of independent, Ber$(x)$-distributed random variables with unknown $x \in [0, 1]$, we prove the existence of an unbiased, $[0, 1]$-valued estimator of the probability $f(x) \in (0, 1)$.
Our estimator is based on so-called debiasing, or randomly truncating a telescopic series of consistent estimators.
Debiased estimators of a probability are not typically constrained to $[0, 1]$, or even bounded, even when all consistent estimators used as inputs are.
We show that constructing the series of consistent estimators from the coefficients of a particular Bernoulli factory yields provable boundedness provided $f \in C^{\rho}[0, 1]$ for $\rho > 3$.
Our result can be thought of as a novel Bernoulli factory with the appealing property that the required number of Ber$(x)$-distributed random variates is independent of their outcomes.
\end{abstract}



\section{Introduction}

Estimating a probability, or generating events with a given probability, are prototypical tasks in probability and computational statistics.
An example is the so-called Bernoulli factory problem: construct a $\text{Ber}( f( x ) )$-distributed random variable (or an $f( x )$-coin) from an almost surely finite number of independent $x$-coins, for a known function $f : [ 0, 1 ] \to [ 0, 1 ]$ but without knowing $x \in [ 0, 1 ]$.
A necessary and sufficient condition for the existence of a Bernoulli factory is for $f$ to be continuous, and either constant or to satisfy
\begin{equation}\label{polynomial_bdd}
\min\{ f( x ), 1 - f( x ) \} \geq \min\{ x, 1 - x \}^n
\end{equation}
for some $n \in \bbN$ \cite{keane/obrien:1994}.
A Bernoulli factory for $f$ can be thought of as a way to generate events with probability $f( x )$, or as a $\{ 0, 1 \}$-valued, unbiased estimator of $f( x )$ produced from a finite number of $\{ 0, 1 \}$-valued, unbiased estimators of $x$.
Other prominent examples include the Metropolis--Hastings algorithm and rejection sampling (see e.g.\ \cite{robert/casella:2004} for an introduction), both of which rely on generating events with a given probability without necessarily evaluating the probability directly.
The approach is prominent in simulation and Bayesian inference under intractable target distributions \cite{lyneetal:2015, stumpf-fetizon/goncalves:2025, vatsetal:2022}, with diffusion models as a prototypical use-case \cite{Beskos06, beskos/roberts:2005, rhee/glynn:2015, sermaidisetal:2013}.

We propose a method for obtaining unbiased estimators $\hat{ f }( x )$ taking values in $[ 0, 1 ]$ for a class of functions $f : [0, 1] \to ( 0, 1 )$, where $x$ is unknown but an infinite sequence of independent $x$-coins is available.
If an unbiased estimator $\hat{ f }( x )$ satisfies $\hat{ f }( x ) \in [ 0, 1 ]$ almost surely, then flipping an $\hat{ f }( x )$-coin yields events with the correct probability \cite[Lemma 2.1]{latuszynskietal:2011}.
Our procedure can thus be viewed both as an unbiased estimator with support $[0, 1]$, and as a novel Bernoulli factory.
From the former point of view our result fits into the family of debiasing methods, which yield unbiased estimators (not necessarily of probabilities) from a consistent sequence of biased estimators \cite{mcleish:2011, rhee/glynn:2015}.
A recent example of this family can be found in \cite{chopinetal:2024}, where $f$ is taken to be analytic and practical methods are obtained by truncating its Taylor series. 
In general unbiased estimators based on resampling can be negative, even when the quantity being estimated and all consistent estimators used as inputs are nonnegative almost surely.
Nonnegative, unbiased estimators of $f( x )$ based only on an almost surely finite number of unbiased estimators of $x$ do not exist in general when $E = \bbR$ \cite[Theorem 2.1]{jacob/thiery:2015}, exist only for appropriately monotone functions when $E = [ a, \infty )$ or $E = ( - \infty, a ]$ for some $a \in \bbR$ \cite[Lemma 3.1]{jacob/thiery:2015}, and exist if and only if
\begin{equation*}
f( x ) \geq \min\{ x - a, b - x \}^n
\end{equation*}
for some $n \in \bbN$ when $E = [ a, b ]$ for $b > a$ \cite[Theorem 3.1]{jacob/thiery:2015}.
Our result extends this body of work to include an almost sure upper bound when $E = [0, 1]$ and $0 < f(x) < 1$ for all $x \in [0,1]$.

From the Bernoulli factory point of view, our method is \textit{non-adaptive}: the number of $x$-coins needed to produce $\hat{ f }( x )$ is independent of their outcomes and can be simulated upfront.
Other non-adaptive Bernoulli factories include that of Keane and O'Brien \cite{keane/obrien:1994} for general functions satisfying \eqref{polynomial_bdd}, and Mendo \cite{mendo:2019} for functions with a convergent power series with nonnegative coefficients.
The factory of Nacu and Peres \cite{holtzetal:2011, nacu/peres:2005} relies on deciding whether to continue flipping $x$-coins based on the outcomes observed thus far.
From an analytical perspective, non-adaptivity was crucial in establishing a recent connection between the existence of Bernoulli factories and duality of certain corresponding pairs of stochastic processes arising in population genetics \cite{koskelaetal:2024}.
It can also be a desirable property for computation, e.g.\ by making it possible to parallelise the simulation task with no need for communication between cores until all $x$-coins have been generated.

\section{Preliminaries}

For $\rho \in \mathbb{N}$ we write $f \in C^\rho[0, 1]$ if $f$ is $\rho$ times continuously differentiable.
For positive $\rho \notin \mathbb{N}$ we let $C^\rho[0, 1]$ stand for the Hölder space that consists of $\lfloor \rho \rfloor$ times continuously differentiable functions whose $\lfloor \rho \rfloor$th derivatives are Hölder continuous with exponent $\rho - \lfloor \rho \rfloor \in (0, 1)$.
That is, $f \in C^\rho[0,1]$ for non-integer $\rho$ if $f$ is $\lfloor \rho \rfloor$ times continuously differentiable and the $\lfloor \rho \rfloor$th derivative $f^{(\lfloor \rho \rfloor)}$ satisfies the Hölder condition
\begin{equation*}
\sup_{x, y \in [0, 1], x \neq y} \frac{ \lvert f^{(\lfloor \rho \rfloor)}(x) - f^{(\lfloor \rho \rfloor)}(y) \rvert }{ \lvert x - y \rvert^{\rho - \lfloor \rho \rfloor} } < \infty.
\end{equation*}
Note that the spaces $C^\rho[0, 1]$ are nested: $C^{\rho_1}[0,1] \subset C^{\rho_2}[0,1]$ for any $0 < \rho_2 < \rho_1$.

Let $x \in [0, 1]$ and let $(X_1, X_2, \ldots )$ be independent $x$-coins.
Let $L$ be a random variable independent of all the coins with $\bbP( L = n ) = \mathds{1}_{ \{ n \geq k \} } n^{ - \lambda } / \zeta( \lambda, k )$, where $k \in \bbN$ and $\lambda > 1$ will be specified later, and
\begin{equation*}
\zeta( \lambda, k ) \coloneqq \sum_{ j = 0 }^{ \infty } ( j + k )^{ - \lambda }
\end{equation*}
is the generalised Riemann zeta function.
For future convenience, we record the following lemma about $\zeta(\lambda, k)$, which follows from \cite[eq.\ (1.3) and Theorem 1.3]{nemes:2017} with $N = 1$ in their notation, and using the fact that the 2nd Bernoulli number is $B_2 = 1/6$.
\begin{lemma}\label{zeta_lemma}
For $\lambda > 1$ and $k > 0$,
\begin{equation*}
\Bigg| \zeta(\lambda, k) - \frac{1}{(\lambda - 1) k^{\lambda - 1}} \Bigg| \leq  \frac{1}{2k^{\lambda}} + \frac{\lambda}{12 k^{1 + \lambda}}.
\end{equation*}
\end{lemma}

Let $S_n \coloneqq X_1 + \ldots + X_n$.
Our main result will rely on \cite[Theorem 8]{holtzetal:2011}, which we slightly rephrase below in a way that will be convenient for our purposes.

\begin{thm}[Theorem 8 of \cite{holtzetal:2011}]\label{holzetal_theorem}
Let $f : [ 0, 1 ] \to ( 0, 1 )$ and let $\rho \in \bbR_+ \setminus \bbN$.
If  $f \in C^{ \rho }[ 0, 1 ]$ then there exist sequences of polynomials
\begin{align}
g_n(x) &= \sum_{ k = 0 }^n \binom{n}{k} a(n, k) x^k (1 - x)^{n - k}, \label{bernstein_lb}\\
h_n(x) &= \sum_{ k = 0 }^n \binom{n}{k} b(n, k) x^k (1 - x)^{n - k}, \notag
\end{align}
such that $0 \leq a( n, k ) \leq b( n, k ) \leq 1$ for all $0 \leq k \leq n < \infty$, 
\begin{align}
a(n, k) &\geq \sum_{ i = 0 }^k \frac{ \binom{ n - m }{ k - i } \binom{ m }{ i } }{ \binom{ n }{ k } } a(m, i) \eqqcolon H_n(m, k), \label{hypergeometric} \\
b(n, k) &\leq \sum_{ i = 0 }^k \frac{ \binom{ n - m }{ k - i } \binom{ m }{ i } }{ \binom{ n }{ k } } b(m, i), \label{hypergeometric_ub}
\end{align}
for any $m \leq n$, and 
\begin{equation}\label{g_h_limits}
\lim_{ n \to \infty } g_n(x) = f(x) = \lim_{ n \to \infty } h_n( x ).
\end{equation}
Moreover,
\begin{equation}\label{exp_diff_bound}
h_n(x) - g_n(x) \leq C n^{ - \rho / 2 }
\end{equation}
for some $C > 0$, uniformly in $x \in [0, 1]$.
\end{thm}
As discussed in~\cite[Section~7]{holtzetal:2011}, the natural extension of Theorem \ref{holzetal_theorem} to integer $\rho$ should involve generalised Lipschitz spaces~\cite[§9 in Chapter~2]{DeVoreLorentz1993}, which are known to be fully characterised by the rate of convergence of polynomial approximation, rather than the more familiar spaces $C^\rho[0,1]$.

\begin{remark}
The statement of \cite[Theorem 8]{holtzetal:2011} gives a certain monotonicity condition on the sequences $\{ g_n \}_{ n \geq 1 }$ and $\{ h_n \}_{ n \geq 1 }$ instead of \eqref{hypergeometric} and \eqref{hypergeometric_ub} (see condition (iv) of their Result 3).
The fact that their monotonicity is equivalent to \eqref{hypergeometric} and \eqref{hypergeometric_ub} is shown on page 97 of \cite{nacu/peres:2005}; see their equation (2) in particular.
\end{remark}

\begin{remark}
Note that the function $f$ in Theorem \ref{holzetal_theorem} satisfies \eqref{polynomial_bdd} by construction since its range is the open set $(0, 1)$.
Theorem \ref{holzetal_theorem} implies the existence of a Bernoulli factory \cite[Theorem 8]{holtzetal:2011}, for which \eqref{polynomial_bdd} is a necessary condition \cite{keane/obrien:1994}.
Thus, the assumption that the range of $f$ is $(0, 1)$ in Theorem \ref{holzetal_theorem} only imposes an additional constraint on $f(0)$ and $f(1)$, which cannot take values in $\{0, 1\}$.
\end{remark}

\section{A debiased Bernoulli factory}

We are now ready to state and prove our main result.

\begin{thm}\label{debiased_bf_theorem}
Let $f : [ 0, 1 ] \to ( 0, 1 )$ and $f \in C^{ \rho }[ 0, 1 ]$ for $\rho > 3$.
Then, for $1 < \lambda < ( \rho - 1 ) / 2$ and any sufficiently large $k \in \mathbb{N}$,
\begin{equation}\label{estimator}
\psi( L, S_L ) \coloneqq H_L( k - 1, S_L ) + \sum_{ n = k }^L \frac{ H_L( n, S_L ) - H_L( n - 1, S_L ) }{ \bbP( L \geq n ) }
\end{equation}
is an unbiased estimator of $f( x )$ for any $x \in [ 0, 1 ]$, and $\psi( L, S_L ) \in [0, 1]$ almost surely.
\end{thm}

\begin{remark}
Our proof relies on controlling the growth of the sum in \eqref{estimator} by applying \eqref{exp_diff_bound} to its individual summands.
Since all summands in \eqref{estimator} are nonnegative, as we will show below, it is clear that $\rho > 2$ is a minimal requirement for finiteness of \eqref{estimator} without further assumptions which make it possible to improve the rate in \eqref{exp_diff_bound}.
We require one more derivative because of the bound in \eqref{binomial_lb}.
Since the latter is tight up to improving constant terms in Stirling's formula, we expect that $\rho > 3$ derivatives are necessary, as well as sufficient, for Theorem \ref{debiased_bf_theorem} to hold.
\end{remark}

\begin{proof}
The proof proceeds in three parts, which consist of establishing that \eqref{estimator} is i) unbiased, ii) nonnegative, and iii) no greater than one.
These three parts prove the theorem by \cite[Lemma 2.1]{latuszynskietal:2011}.
Each part will make use of the conclusions of Theorem \ref{holzetal_theorem}, whose hypotheses are satisfied by assumption, and by nestedness of H\"older spaces: $f \in C^{ \rho }[ 0, 1 ] \Rightarrow f \in C^{ \eta }[ 0, 1 ]$ for any $\eta < \rho$, so that we may assume $\rho \notin \bbN$.

For lack of bias, it suffices that for $n \leq L$,
\begin{align}
\bbE[ H_L( n, S_L ) | L ] &= \sum_{ s = 0 }^L \sum_{ i = 0 }^s \binom{ L - n }{ s - i } \binom{ n }{ i } a(n, i) x^s (1 - x)^{ L - s } \notag \\
&= \sum_{ i = 0 }^n a(n, i) \binom{ n }{ i } x^i (1 - x)^{ n - i } \sum_{ s = i }^{ L - n + i } \binom{ L - n }{ s - i } x^{ s - i } (1 - x)^{ L - s - n + i } \notag \\
&= \sum_{ i = 0 }^n a(n, i) \binom{ n }{ i } x^i (1 - x)^{ n - i } \sum_{ s = 0 }^{ L - n } \binom{ L - n }{ s } x^s (1 - x)^{ L - n - s } \notag \\
&= \sum_{ i = 0 }^n a(n, i) \binom{ n }{ i } x^i (1 - x)^{ n - i } = g_n( x ), \label{lack_of_bias}
\end{align}
which amounts to saying that the mean of $n$ $x$-coins coincides with that of a size-$n$ subsample picked uniformly from $L$ $x$-coins.
The result follows as in \cite[Theorem 1]{rhee/glynn:2015}:
\begin{align*}
\bbE[ \psi( L, S_L ) ] &= \bbE\Bigg[ H_L( k - 1, S_L ) + \sum_{ n = k }^L \frac{ H_L( n, S_L ) - H_L( n - 1, S_L ) }{ \bbP( L \geq n ) } \Bigg] \\
&= \bbE[ H_L( k - 1, S_L ) ] + \bbE\Bigg[ \sum_{ n = k }^{ \infty } \frac{ [ H_L( n, S_L ) - H_L( n - 1, S_L ) ] \mathds{ 1 }\{ L \geq n \}  }{ \bbP( L \geq n ) } \Bigg] \\
&= \bbE[ H_L( k - 1, S_L ) ] + \sum_{ n = k }^{ \infty } \bbE\Bigg[ \frac{ \mathds{ 1 }\{ L \geq n \}  }{ \bbP( L \geq n ) } \bbE[ H_L( n, S_L ) - H_L( n - 1, S_L ) | L ]  \Bigg] \\
&= g_{ k - 1 }( x ) + \sum_{ n = k }^{ \infty } [ g_n( x ) - g_{ n - 1 }( x ) ] = f(x),
\end{align*}
where the final equality uses \eqref{g_h_limits}.

Nonnegativity follows from \eqref{hypergeometric}: for $m \leq n$,
\begin{align}
&H_n( m, k ) - H_n( m - 1, k ) \notag \\
&= \sum_{ i = 0 }^k \Bigg( \frac{ \binom{ n - m }{ k - i } \binom{ m }{ i } }{ \binom{ n }{ k } } a( m, i ) - \frac{ \binom{ n - m + 1 }{ k - i } \binom{ m - 1 }{ i } }{ \binom{ n }{ k } } a( m - 1 , i ) \Bigg) \nonumber \\
&\geq \sum_{ i = 0 }^k \Bigg( \frac{ \binom{ n - m }{ k - i } \binom{ m }{ i } }{ \binom{ n }{ k } } \Big[ \frac{ m - i }{ m } a( m - 1, i ) + \frac{ i }{ m } a( m - 1 , i - 1 )\Big]  - \frac{ \binom{ n - m + 1 }{ k - i } \binom{ m - 1 }{ i } }{ \binom{ n }{ k } } a( m - 1 , i ) \Bigg) \nonumber \\
&= \sum_{ i = 0 }^k \Bigg( \frac{ \binom{ n - m }{ k - i } \binom{ m }{ i } }{ \binom{ n }{ k } } \frac{ m - i }{ m } + \frac{ \binom{ n - m }{ k - i - 1 } \binom{ m }{ i + 1 } }{ \binom{ n }{ k } } \frac{ i + 1 }{ m } - \frac{ \binom{ n - m + 1 }{ k - i } \binom{ m - 1 }{ i } }{ \binom{ n }{ k } } \Bigg) a( m - 1, i ). \label{non_neg}
\end{align}
Now
\begin{align}
&\binom{ n - m }{ k - i } \binom{ m }{ i } \frac{ m - i }{ m } + \binom{ n - m }{ k - i - 1 } \binom{ m }{ i + 1 } \frac{ i + 1 }{ m } - \binom{ n - m + 1 }{ k - i } \binom{ m - 1 }{ i } \nonumber \\
&= \binom{ n - m + 1 }{ k - i } \binom{ m - 1 }{ i } \Bigg( \frac{ n - m - k + i + 1 }{ n - m + 1 } + \frac{ k - i }{ n - m + 1 } - 1 \Bigg) = 0. \label{combinatorics}
\end{align}
Hence, the right-hand side of \eqref{non_neg} vanishes, which makes the left-hand side, and thus every summand in \eqref{estimator}, nonnegative.

To obtain an upper bound, we begin similarly:
\begin{align}
&H_n( m, k ) - H_n( m - 1, k ) \notag \\
&= \sum_{ i = 0 }^k \Bigg( \frac{ \binom{ n - m }{ k - i } \binom{ m }{ i } }{ \binom{ n }{ k } } a( m, i ) - \frac{ \binom{ n - m + 1 }{ k - i } \binom{ m - 1 }{ i } }{ \binom{ n }{ k } } a( m - 1 , i ) \Bigg) \nonumber \\
&\leq \sum_{ i = 0 }^k \Bigg( \frac{ \binom{ n - m }{ k - i } \binom{ m }{ i } }{ \binom{ n }{ k } } \Big[ \frac{ m - i }{ m } b( m - 1, i ) + \frac{ i }{ m } b( m - 1 , i - 1 )\Big] - \frac{ \binom{ n - m + 1 }{ k - i } \binom{ m - 1 }{ i } }{ \binom{ n }{ k } } a( m - 1 , i ) \Bigg) \nonumber \\
&= \sum_{ i = 0 }^k \Bigg[\Bigg( \frac{ \binom{ n - m }{ k - i } \binom{ m }{ i } }{ \binom{ n }{ k } } \frac{ m - i }{ m } + \frac{ \binom{ n - m }{ k - i - 1 } \binom{ m }{ i + 1 } }{ \binom{ n }{ k } } \frac{ i + 1 }{ m } \Bigg) b( m - 1, i ) - \frac{ \binom{ n - m + 1 }{ k - i } \binom{ m - 1 }{ i } }{ \binom{ n }{ k } } a( m - 1, i ) \Bigg] \notag \\
&= \sum_{ i = 0 }^k \frac{ \binom{ n - m + 1 }{ k - i } \binom{ m - 1 }{ i } }{ \binom{ n }{ k } } [ b( m - 1, i ) - a( m - 1, i ) ], \label{b_minus_a}
\end{align}
where the inequality follows from $a(n, k) \leq b(n, k)$ and \eqref{hypergeometric_ub}, and the last line from \eqref{combinatorics}.
To bound \eqref{b_minus_a} from above we need to control $b(n, k) - a(n, k)$, which we do next.

By uniformity of \eqref{exp_diff_bound} we have the trivial bounds
\begin{align}
b(n,0) - a(n,0) &= h_n(0) - g_n(0) \leq C n^{ - \rho / 2 }, \label{trivial_0_bound} \\
b(n,n) - a(n,n) &= h_n(1) - g_n(1) \leq C n^{ - \rho / 2 }. \label{trivial_1_bound}
\end{align}
For $k \in \{1, \ldots, n - 1\}$, we also have
\begin{align*}
b(n,k) - a(n,k) &= \frac{ \binom{ n }{ k } ( k / n )^k ( 1 - k / n )^{ n - k } [ b( n, k ) - a( n, k) ] }{ \binom{ n }{ k } ( k / n )^k ( 1 - k / n )^{ n - k } } \\
&\leq \sum_{j = 0}^n \frac{ \binom{ n }{ j } ( k / n )^j ( 1 - k / n )^{ n - j } [ b( n, j ) - a( n, j ) ] }{ \binom{ n }{ k } ( k / n )^k ( 1 - k / n )^{ n - k } } \\
&= \frac{ h_n( k / n ) - g_n( k / n ) }{ \binom{ n }{ k } ( k / n )^k ( 1 - k / n )^{ n - k } } \leq \frac{ C n^{ -\rho / 2 } }{ \binom{ n }{ k } ( k / n )^k ( 1 - k / n )^{ n - k } },
\end{align*}
where the first step multiplies and divides by the same constant and the inequality adds nonnegative terms.
We bound the denominator on the right-hand side from below using the non-asymptotic Stirling approximations \cite{robbins:1955},
\begin{equation*}
\sqrt{ 2 \pi } n^{ n + 1/2 } e^{ -n } \leq  n! \leq \sqrt{ 2 \pi } n^{ n + 1/2 } e^{ 1 - n },
\end{equation*}
to obtain
\begin{equation}\label{binomial_lb}
b(n,k) - a(n,k) \leq \frac{ C n^{ -\rho / 2} }{ \sqrt{ \frac{ n }{ 2 \pi k ( n - k ) } } e^{-2} } \leq C \sqrt{ \frac{ \pi }{ 2 } } e^2 n^{ (1 - \rho) / 2 },
\end{equation}
where the last step follows by evaluating the denominator in the middle step at the unique global minimum $k = n / 2$.
Since the right-hand side of \eqref{binomial_lb} is greater than that of \eqref{trivial_0_bound} or \eqref{trivial_1_bound}, plugging it into \eqref{b_minus_a} yields
\begin{equation*}
H_n( m, k ) - H_n( m - 1, k ) \leq C \sqrt{ \frac{ \pi }{ 2 } } e^2 ( m - 1 )^{ ( 1 - \rho ) / 2 },
\end{equation*}
whereupon substituting in the survival function $\bbP( L \geq n ) = \zeta( \lambda, n ) / \zeta( \lambda, k )$ and using
\begin{equation*}
(n - 1)^{(1 - \rho) / 2} = \frac{n^{(1 - \rho) / 2}}{(1 - 1 / n)^{(\rho - 1)/2}}
\end{equation*}
results in the bound 
\begin{align}
\psi( L, S_L ) &= H_L( k - 1, S_L ) + \sum_{ n = k }^L \frac{ H_L( n, S_L ) - H_L( n - 1, S_L ) }{ \mathbb{P}( L \geq n ) } \notag \\
&\leq H_L( k - 1, S_L ) + C \sqrt{ \frac{ \pi }{ 2 } } e^2 \frac{\zeta( \lambda, k )}{(1 - 1 / k)^{(\rho - 1)/2}} \sum_{ n = k }^{ \infty } \frac{ n^{ ( 1 - \rho ) / 2 } }{ \zeta( \lambda, n ) }. \label{hurwitz_bound}
\end{align}
Since $\zeta( \lambda, n ) = \Theta( n^{ - \lambda + 1 } )$ by Lemma \ref{zeta_lemma}, the sum on the right-hand side of \eqref{hurwitz_bound} is finite if 
\begin{equation*}
\frac{ 1 }{ 2 } - \frac{ \rho }{ 2 } + \lambda - 1 < -1 \Rightarrow \lambda < \frac{ \rho }{ 2 } - \frac{ 1 }{ 2 },
\end{equation*}
which can be satisfied by $\lambda > 1$ since $\rho > 3$ by assumption.
For such a $\lambda$, Lemma \ref{zeta_lemma} shows that \eqref{hurwitz_bound} satisfies the bound
\begin{align*}
&C \sqrt{ \frac{ \pi }{ 2 } } e^2 \frac{\zeta( \lambda, k )}{(1 - 1 / k)^{(\rho - 1)/2}} \sum_{ n = k }^{ \infty } \frac{ n^{ ( 1 - \rho ) / 2 } }{ \zeta( \lambda, n ) } \\
&\leq \frac{C e^2 \sqrt{\pi / 2}}{(1 - 1 / k)^{(\rho - 1)/2}} \Big(1 + \frac{\lambda  - 1}{2 k}\Big(1 + \frac{\lambda}{6 k}\Big)\Big) \frac{1}{k^{\lambda - 1}} \sum_{ n = k }^{ \infty } \frac{n^{\lambda - (1 + \rho) / 2}}{1 - \frac{\lambda  - 1}{2 n}(1 + \frac{\lambda}{6 n})} \\
&\leq \frac{C e^2 \sqrt{\pi / 2}}{(1 - 1 / k)^{(\rho - 1)/2}} \frac{2k + (\lambda  - 1)(1 + \frac{\lambda}{6 k})}{2k - (\lambda  - 1)(1 + \frac{\lambda}{6 k})} \frac{\zeta((1 + \rho) / 2 - \lambda, k)}{k^{\lambda - 1}} \\
&\leq \frac{C e^2 \sqrt{\pi} (2k + (\lambda  - 1)(1 + \frac{\lambda}{6 k}))(2k + ((\rho - 1) / 2 - \lambda)(1 + \frac{(\rho + 1) / 2 - \lambda}{6 k}))}{\sqrt{2} (1 - 1 / k)^{(\rho - 1)/2} (2k - (\lambda  - 1)(1 + \frac{\lambda}{6 k})) 2 k ((\rho - 1) / 2 - \lambda)} \frac{1}{k^{(\rho - 3) / 2}}.
\end{align*}
To conclude the proof, it remains to bound $H_L( k - 1, S_L )$ sufficiently far away from 1 with probability 1 so that \eqref{hurwitz_bound} can be made less than one.
We will do this by noting that by \eqref{hypergeometric},
\begin{equation*}
H_L(k - 1, S_L) \leq a(L, S_L) \leq b(L, S_L),
\end{equation*}
and showing that there exists $\delta > 0$ such that $\max_{k \in \{0, \ldots, n\}} b(n, k) \leq 1 - \delta$ for any sufficiently large $n$.

Consider the upper polynomial approximant $h_n(x)$ defined in \cite[eq.\ (48)]{holtzetal:2011}, which satisfies all conditions of Theorem \ref{holzetal_theorem} by construction:
\begin{equation*}
h_n(x) = f_n(x) + D \Big(\frac{\theta}{n^{\rho}} + \Big(\frac{x(1-x)}{n}\Big)^{\rho / 2} \Big),
\end{equation*}
where $D > 0$ and $\theta > 0$ are constants and $f_n(x)$ is defined by iteration \cite[page 352]{holtzetal:2011}.
By \cite[Lemma 25]{holtzetal:2011}, there is a constant $C_{\theta} > 0$ such that
\begin{equation}\label{theta_and_D}
D \Big(\frac{\theta}{n^{\rho}} + \Big(\frac{x(1-x)}{n}\Big)^{\rho / 2} \Big) \leq \frac{C_{\theta} D}{n^{\rho / 2}},
\end{equation}
and by the paragraph labelled \textbf{Step 4.}\ in \cite[pages 354--355]{holtzetal:2011}, the Bernstein coefficients of $f_n$ are bounded above by $1 - \varepsilon / 2$ whenever $f \in [\varepsilon, 1 - \varepsilon]$.
By continuity, such an $\varepsilon > 0$ exists for our $f : [0, 1] \to (0, 1)$.
By linearity of the Bernstein operator, we thus have that the Bernstein coefficients of $h_n(x)$ satisfy
\begin{equation*}
b(n, k) \leq 1 - \frac{\varepsilon}{2} + \frac{C_{\theta} D}{n^{\rho / 2}}
\end{equation*}
for any  $k \in \{0, \ldots, n\}$.
Thus
\begin{align*}
H_L( k - 1, S_L ) = \sum_{ i = 0 }^{ S_L } \frac{ \binom{ L - k + 1 }{ S_L - i } \binom{ k - 1 }{ i } }{ \binom{ L }{ S_L } } a( k - 1, i ) \leq 1 - \frac{\varepsilon}{2} + \frac{C_{\theta} D}{(1-1/k)^{\rho/2}}\frac{1}{k^{\rho / 2}},
\end{align*}
which completes the proof since 
\begin{align}
&\psi(L, S_L) \notag \\
&\leq 1 - \frac{\varepsilon}{2} + \frac{C_{\theta} D}{(1-1/k)^{\rho/2}}\frac{1}{k^{\rho / 2}} \notag\\
&\phantom{\leq 1 -} + \frac{C e^2 \sqrt{\pi} (2k + (\lambda  - 1)(1 + \frac{\lambda}{6 k}))(2k + ((\rho - 1) / 2 - \lambda)(1 + \frac{(\rho + 1) / 2 - \lambda}{6 k}))}{\sqrt{2} (1 - 1 / k)^{(\rho - 1)/2} (2k - (\lambda  - 1)(1 + \frac{\lambda}{6 k})) 2 k ((\rho - 1) / 2 - \lambda)} \frac{1}{k^{(\rho - 3) / 2}} \notag \\
&< 1 \label{upper_bound_decay}
\end{align}
for any sufficiently large $k$.
\end{proof}

By discarding those assumptions from Theorem \ref{debiased_bf_theorem} which were only used in the proof of the upper bound of 1, we obtain the following corollary.
\begin{cor}\label{nonneg_corollary}
Let $f : [0, 1] \to (0, 1)$ and suppose there are coefficients $0 \leq a(n, k) \leq 1$ such that the functions $g_n(x)$ defined in \eqref{bernstein_lb} satisfy $\lim_{n \to \infty} g_n(x) = f(x)$.
Then, for any $\lambda > 1$ and any $k \in \bbN$, the estimator $\psi( L, S_L )$ defined in \eqref{estimator} is unbiased and non-negative almost surely.
\end{cor}
An analogous calculation shows that a version of \eqref{estimator} obtained by defining $H_n(m, k)$ as the right-hand side of \eqref{hypergeometric_ub} rather than of \eqref{hypergeometric} is unbiased and bounded above by one.
Non-negative estimators of probabilities find applications in pseudo-marginal MCMC \cite{Andrieu09} and can be challenging to obtain for so-called doubly-intractable problems \cite{lyneetal:2015}.
We are not aware of any particular uses for estimators bounded only from above.

Because $\bbP( \psi( L, S_L ) \in [ 0, 1 ] ) = 1$ and the maximal variance of a $[0, 1]$-supported random variable is $1/4$, we have a further corollary.
\begin{cor}
Under the assumptions of Theorem \ref{debiased_bf_theorem}, $\text{Var}( \psi( L, S_L ) ) \leq 1/4$.
\end{cor}

It is unclear whether $\rho > 3$ is necessary for \eqref{estimator} to have finite variance.
It is straightforward to verify that $2 < \rho$ is sufficient for the variance formula in \cite[Eq.\ (6)]{rhee/glynn:2015} to converge:
\begin{align*}
&\sum_{n = 1}^{\infty} \frac{\bbE[(H_L(n - 1, S_L) - f(x))^2 - (H_L(n - 1, S_L) - f(x))^2]}{\bbP(L \geq n)} \\
&\leq \sum_{n = 1}^{\infty} \frac{\bbE[H_L(n - 1, S_L)^2 - H_L(n - 1, S_L)^2] + 2 C n^{-\rho / 2}}{\bbP(L \geq n)},
\end{align*}
where we have expanded both squares, and used \eqref{lack_of_bias} as well as \eqref{exp_diff_bound}.
By essentially the same calculation leading to \eqref{non_neg}, we further have that $H_L(n - 1, S_L)^2 - H_L(n - 1, S_L)^2 \leq 0$, so that, for a constant $\tilde{C} > 0$ independent of $n$, 
\begin{equation*}
\sum_{n = 1}^{\infty} \frac{\bbE[(H_L(n - 1, S_L) - f(x))^2 - (H_L(n - 1, S_L) - f(x))^2]}{\bbP(L \geq n)} \leq \tilde{C} \sum_{n = 1}^{\infty} n^{-\rho / 2 + \lambda - 1},
\end{equation*}
which is finite provided $2 < 2 \lambda < \rho$.
However, the derivation of \cite[Eq.\ (6)]{rhee/glynn:2015} assumes square-integrability of $\psi(L, S_L)$, for which decay of the form
\begin{equation}\label{finite_variance}
\lim_{n \to 0} \frac{\bbE[(H_L(n, S_L) - f(x))^2]}{\bbP(L \geq n)} = 0
\end{equation}
seems necessary.
We have not been able to verify whether $\rho > 2$ is sufficient for \eqref{finite_variance}, and hence whether it also suffices for finite variance of the non-negative estimators in Corollary \ref{nonneg_corollary}.

\section{Implementability}\label{implementation}

The drawback of the Nacu--Peres factory is that it is often difficult to construct the sequences $a(n, k)$ and $b(n, k)$ needed in Theorem \ref{debiased_bf_theorem} for a given function $f$.
A formula is provided in \cite[eq.\ (48)]{holtzetal:2011}, but it is based on an iterative procedure \cite[page 352]{holtzetal:2011}.
An example of a similar iteration was analysed in \cite{latuszynskietal:2011}, and would have required evaluation of $2^{2^{26}}$ coefficients.
The authors of \cite{latuszynskietal:2011} presented a generalisation in which construction of a Bernoulli factory only required analogues of \eqref{hypergeometric} and \eqref{hypergeometric_ub} to hold for certain conditional expectations \cite[eq.\ (6) and (7)]{latuszynskietal:2011}, resulting in much more practical algorithms.
In our setting the same generalisation would guarantee only that the same conditional expectations of $\psi( L, S_L )$ lie in $[0, 1]$, which is insufficient for our purposes.

Even when the coefficients $a(n, k)$ and $b(n, k)$ are available, our method requires the specification of the initial value $k$ in \eqref{estimator}, corresponding to the least element in the support of $L$.
The least element has to be large enough to ensure \eqref{upper_bound_decay} holds, which depends on five constants: $\varepsilon$, which is straightforward to determine from $f$; $\theta$, $C_{\theta}$ and $D$ from \eqref{theta_and_D}; and $C$ from \eqref{exp_diff_bound} which can be taken to be $C = 2 C_{\theta} D$ \cite[eq.\ (48)]{holtzetal:2011}.
The constants $C_{\theta}$ and $D$ are complicated, and require tracking various other constants (some defined only implicitly) through the calculation labelled \textbf{Step 5.}\ in \cite[pages 355--357]{holtzetal:2011}.
Affine functions $f = a + b x$ satisfying \eqref{polynomial_bdd} constitute an exception since Bernstein polynomials reproduce them exactly, so that it is always possible to take $k = 1$, $a(n, k) = b(n, k) = f(k / n)$, and $\bbP(L = 1) = 1$ so that the estimator reduces to $\psi(1, S_1) = a + \mathds{1}\{S_1 = 1\} b$.
In comparison, the non-adaptive Keane--O'Brien estimator \cite{keane/obrien:1994} is not straightforward to compute even for this special case because it relies on an iterative construction of a sequence of sets of binary strings whose definition is implicit.

On the other hand, the non-negative estimators from Corollary \ref{nonneg_corollary} (resp.\ their counterparts bounded only from above discussed immediately afterwards) are readily implementable as soon as the coefficients $a(n, k)$ (resp.\ $b(n, k)$) can be computed.
Since Corollary \ref{nonneg_corollary} does not require a specific rate of decay of $f(x) - g_n(x)$ as in \eqref{exp_diff_bound}, it suffices to take e.g.\ $a(n, k) = f(k / n) - (1 + \sqrt{2}) M / \sqrt{n}$ when $f$ is $M$-Lipschitz, or $a(n, k) = f(k / n) - \| f'' \|_{\infty} / (2n)$ if $f \in C^2[0, 1]$ \cite[Proposition 10]{nacu/peres:2005}.
In both cases, the least element of the support of $L$ can be taken as any value $n_0$ which ensures $a(n, k) \geq 0$ for all $n \geq n_0$.
However, neither of these examples has sufficient smoothness to guarantee finite variance, irrespective of whether $\rho > 2$ or $\rho > 3$ is required.

\section{Open problems}

The number of $x$-coins needed to output an $f(x)$-coin is a natural efficiency criterion for a Bernoulli factory.
Hence, it is desirable to choose $\lambda$ close to the upper limit of $( \rho - 1 ) / 2$ in Theorem \ref{debiased_bf_theorem} to make the tail of the random number of coins as light as possible.
We conjecture that exponential tail decay can be obtained for an appropriate class of $f \in C^{ \infty }[ 0, 1 ]$-functions.
The necessary ingredient for adapting the proof of Theorem \ref{debiased_bf_theorem} would be a bound of the form
\begin{equation} \label{exp_diff_bound_2}
h_n( x ) - g_n( x ) < C \gamma^n
\end{equation}
for $C > 0$ and $ \gamma \in ( 0, 1 )$ uniformly in $x \in [ 0, 1 ]$.
Such a bound is available for any closed $E \subset (0, 1)$ \cite[Theorem 2]{nacu/peres:2005}.
However, excluding the $x \in \{0, 1\}$ boundaries rules out the bounds in \eqref{trivial_0_bound} and \eqref{trivial_1_bound}, which are required because all $x$-coins can take identical values even when $x \in (0, 1)$.
The construction of polynomials $h_n$ and $g_n$ that satisfy \eqref{exp_diff_bound} in \cite{holtzetal:2011} depends on $\lfloor \rho \rfloor$, as is seen in their equations (47)--(48).
Therefore an exponential bound akin to \eqref{exp_diff_bound_2} requires a different construction.

It would also be of interest to extend the results of Theorem \ref{debiased_bf_theorem} to functions $f : [0, 1] \to [0, 1]$, or $f : (0, 1) \to (0, 1)$.
The former would require adapting \cite[Theorem 8]{holtzetal:2011} to handle boundary points $f(0), f(1) \in \{0, 1\}$, which appears difficult since the construction of their approximating polynomials in their Section 5 relies on $0 < \min_{x \in [0, 1]}\{f(x)\} \leq \max_{x \in [0, 1]}\{f(x)\} < 1$.
The difficulty with the latter lies in establishing \eqref{trivial_0_bound} and \eqref{trivial_1_bound}, which amounts to controlling $h_n(x) - g_n(x)$ at $x \in \{0, 1\}$ even when the domain of $f$ is $(0, 1)$.
This can be done trivially when $f$ can be continuously extended to a function $f : [0, 1] \to (0, 1)$.
When that isn't possible, handling $f : (0, 1) \to (0, 1)$ seems to be as difficult as $f : [0, 1] \to [0, 1]$.

\section*{Acknowledgements and data sharing}

We thank Francesca Crucinio for initiating the discussion on the variance of $\psi( L, S_L )$.
This work began while JK was supported by EPSRC research grant EP/V049208/1 at the University of Warwick.
TK was supported by the Research Council of Finland decisions 338567, 359183 and 368086.
K{\L} was supported by the Royal Society through the Royal Society University Research Fellowship.
Data sharing is not applicable to this article as no new data were generated or analysed.

\bibliographystyle{alpha}
\bibliography{bibliography}  

\end{document}